\documentclass[12pt]{amsart}

\usepackage{graphicx}
\usepackage{amsmath, comment}
\usepackage{amscd}
\usepackage{amsfonts}
\usepackage{amssymb}
\usepackage{color}
\usepackage{fullpage}
\usepackage{mathtools}
\usepackage[hypertexnames=false, colorlinks, citecolor=red,linkcolor=blue, urlcolor=black]{hyperref}
\usepackage{tikz-cd}
\usepackage[capitalise]{cleveref}
\usepackage{enumitem}

\numberwithin{equation}{section}

\newtheorem*{theorem*}{Theorem}
\newtheorem{thm}{Theorem}[section]

\theoremstyle{plain}
\newtheorem{lem}[thm]{Lemma}

\newtheorem{rema}[thm]{Remark}

\makeatletter
\newcommand*{\rom}[1]{\expandafter\@slowromancap\romannumeral #1@}

\theoremstyle{definition}

\begin{document}
\author{Amin Soofiani}   
\address{Department of Mathematics, University of British Columbia, Vancouver, B.C., V6T
1Z2 Canada} 
\email{soofiani@math.ubc.ca}

\title{AN EXPLICIT FORMULA FOR LARMOUR'S DECOMPOSITION
OF HERMITIAN FORMS}

\begin{abstract}
   Let $K$ be a complete discretely valued field whose residue field has characteristic different from $2$. Let $(D,\sigma)$ be a $K-$division algebra with involution of the first kind, and $h$ be a $K-$anisotropic $\epsilon$-hermitian form over $(D,\sigma)$. By a theorem due to Larmour, there is a decomposition $h=h_0 \perp h_1$ such that the elements in a diagonalization of $h_0$ are units, the elements in a diagonalization of $h_1$ are uniformizers, and $h_0$, $h_1$ are determined uniquely up to $K-$isometry. In this paper, we give an explicit description of the elements in the diagonalization of $h_0$ and $h_1$ in the case of quaternion algebras. Then we will derive explicit formulas for Larmour's isomorphism of Witt groups.
\end{abstract} 

\maketitle

\setcounter{tocdepth}{1}
\tableofcontents

\section{Introduction}\label{51}

 let $K$  be a complete discretely valued field with valuation $\nu_K: K \backslash \{0\} \longrightarrow \mathbb{Z}$, uniformizer $\pi$, and residue field $k$ with $char \ k\ne 2$. Let $\mathcal{O}_K$ be the ring of integers of $K$ and $m_K$ be its maximal ideal.

 If a quadratic form over $K$ has a diagonalization with elements in ${\mathcal{O}_K}^*$, then we call it an unramified quadratic form. To any unramified quadratic form $q=<u_1, u_2, \dots, u_r>$ over $K$, one can associate the quadratic form $<\overline{u_1}, \overline{u_2}, \dots, \overline{u_r}>$ over $k$, which we denote by $\overline{q}$.

 Up to isometry, any quadratic form $q$ over $K$ can be decomposed as $q=q_0 \perp \pi q_1$, where $q_0$ and $q_1$ are unramified quadratic forms. Springer's decomposition theorem for quadratic forms (see \cite[Chapter~6, Section~1]{5}) states that if $q$ is $K$-anisotropic, then $q_0$ and $q_1$ are $K$-anisotropic, and the quadratic forms $q_0$ and $\pi q_1$ are determined uniquely up to $K-$isometry. The quadratic forms $\overline{q_0}$ and $\overline{q_1}$ are called the first and second residue forms of $q$, respectively.

\begin{thm}\label{53}
    \textbf{(Springer)} The following map is a group isomorphism:

$$\partial= (\partial_0, \partial_1): W(K) \longrightarrow  W(k) \oplus  W(k)$$
$$\   \  \   [q]=[q_0 \perp \pi q_1]   \     \ \mapsto  \  \   \   {([\overline{q_0}], [\overline{q_1}]}),$$where $W(K)$ and $W(k)$ denote the Witt rings of $K$ and $k$, respectively.

\begin{proof}
    See \cite[Chapter ~6, Theorem ~1.4]{5}.
\end{proof}
  
\end{thm}
We remark that $\partial_1$ depends on the choice of $\pi$, but $\partial_0$ does not.

Now let $(D,\sigma)$ be a $K-$division algebra with involution of the first kind. 
We use the notation "$\epsilon-$hermitian" forms over $(D,\sigma)$ for $\epsilon \in \{ 1, -1\}$ (see \cite[Chapter ~1, Section ~2, Page ~8]{3}). For $\epsilon=1$ (respectively $\epsilon=-1$), the forms are also called hermitian (respectively, skew-hermitian). All $\epsilon-$hermitian forms over $(D,\sigma)$ are diagonalizable by $\epsilon-$symmetric elements (see \cite[Chapter 1, Section 6, 6.2.1-6.2.3]{3}).  

In \cite{6}, Larmour generalized Springer's theorem to the case of $\epsilon-$hermitian forms over $(D,\sigma)$ over Henselian fields.

\begin{thm}
    \label{25}

(\textbf{Larmour}) The map 

$$\partial= (\partial_0, \partial_1): W^{\epsilon}(D,\sigma) \longrightarrow  W^{\epsilon}(\overline{D},\overline{\sigma}) \oplus  W^1(\overline{D},\overline{\sigma_{\pi'}})$$
$$\   \  \ \ \   \  \    \   \    \   \  [h]   \     \ \mapsto  \  \   \   {([\overline{h_0}], [\overline{h_1}]})$$

is an isomorphism of groups.
\begin{proof}
    See \cite[Theorem 3.7]{6}. 
\end{proof}
\end{thm}

The residue algebra $\overline{D}$, twisted involution $\sigma_{\pi'}$ on $D$, and residue involutions $\overline{\sigma}$ and $\overline{\sigma_{\pi'}}$ will be recalled in details in Section \ref{200}.

In this paper, we focus on Larmour's theorem in the case when $D$ is a quaternion division algebra. In order to compute the image of the $K-$isometry classes of $\epsilon-$hermitian forms under the isomorphism $\partial$, one just needs to compute the image of one dimensional forms. Let $h=<u>$ be a one dimensional $\epsilon-$hermitian form over $(D,\sigma)$, where $u\in Sym^{\epsilon}(D,\sigma)$. There are  pure quaternions $x,y\in D$ such that  $x^2=a$, $y^2=b$, $D=(\frac{a,b}{K})$. Also put $z:=xy=-yx$. The set $\{1,x,y,z\}$ is a $K-$basis of $D$. Hence there are elements $\delta, \alpha, \beta, \gamma \in K$ such that $u=\delta + \alpha x + \beta y + \gamma z$. In Section \ref{20} (Lemma \ref{47}), we will show that up to $K-$isometry, the elements  $\delta, \alpha, \beta, \gamma$ can be chosen from the ring of integers $\mathcal{O}_K$. Moreover, depending on the ramification of $D$, the type of $\sigma$, and the ramification of the element $\pi'\in D$, we  explicitly compute the forms $\overline{h_0}$ and $\overline{h_1}$. We show that there are ten cases that need to be considered separately. In each case we explicitly explain the isomorphism $\partial$.

In studying semisimple linear algebraic groups, one often needs to use their classification as certain groups attached to $\epsilon-$hermitian forms defined over central simple algebras with involution. An important open problem in the theory of linear algebraic groups is the norm principle (See \cite[Section 2]{1}).
In \cite{2}, it was shown that the norm principle for the extended Clifford groups of quadratic forms over a complete discretely valued field $K$ can be reduced to the norm principle for the extended Clifford groups of quadratic forms over the residue field $k$. In \cite{7}, the explicit formulas of this paper are used  to generalize the main result in \cite{2} to the case of the extended Clifford group of skew-hermitian forms defined over quaternion algebras over complete discretely valued fields. 

This paper contains 5 sections, including this introduction. In Section \ref{200}, we recall the extension of valuations from fields to division algebras, which is needed to state Larmour's theorem. Section \ref{22} is devoted to the case of quaternion algebras. In order to study Larmour's isomorphism explicitly, we consider several cases, each case depending on the ramification of the algebra and the type of the involution. In each case, we compute the invariants which are needed in later sections. In Section \ref{20}, we will show that up to isometry, the diagonal entries of forms can be chosen from certain subspaces of the division algebra. Using results in Sections \ref{22} and \ref{20}, we derive explicit formulas for Larmour's isomorphism in Section \ref{52}.

\section{Extension of valuations to division algebras and Larmour's Theorem}\label{200}

Let $K, \nu, \pi, k, \mathcal{O},m$ be defined as in Section \ref{51}. In this section we assume that $(D, \sigma)$ a central $K-$division algebra with involution of the first kind, and $n=deg \ D$. 

In order to study an analog of Springer's decompsotion for $\epsilon-$hermitian forms over division algebras defined over complete discretely valued fields, one should first consider the extensions of valuations from fields to division algebras. Let $\nu_K: K \backslash \{0\} \longrightarrow \mathbb{Z}$ denote the valuation on $K$. By \cite{9}, $\nu_K$ extends uniquely to the following valuation on $D$: 

$$\nu_{D}: D \backslash \{0\}   \longrightarrow \mathbb{Q}$$

$$ \   \  \   \   \  \   \   \   \   \  \  \   \  \  \   \   \  \   \   \   \   \  \  \   \  \  \  \  \  \     x \mapsto \frac{1}{n} \  \nu_K (Nrd (x)).$$

 Let $\mathcal{O}_D:= \{d \in D  \ | \ \nu_{D}(d) \geq 0\}$ be the ring of integers of $D$, $m_D:=\{d \in D \  | \  \nu_{D}(d) > 0\}$ its maximal ideal, and $\Gamma_D$ the image of $\nu_D$. 
 
 \[\begin{CD}
 m_K @>>> \mathcal{O}_K @>>> K \\
 @VVV @VVV @VVV     \\
 m_D @>>>  \mathcal{O}_D @>>> D 
\end{CD}
\]\  

All arrows are inclusion in the diagram above.

 Let $j=\frac{1}{min (\Gamma_D \cap  \mathbb{Q}^+)}$, which measures the ramification of $\nu$ on $D$. One extreme case is when $j=1$, in which case we call $D$ unramified. Choose some $\pi'\in D^*$ such that $\nu_D(\pi')=\frac{1}{j}$. The element $\pi'$ is in fact a uniformizer of $D$. Note that $j$ and $\pi'$ do not depend on the involution $\sigma$.

 Fix $\epsilon \in \{1, -1\}$ and choose  an $\epsilon-$symmetric element $\pi'' \in m_D$ with minimal value, i.e. for all $\epsilon-$symmetric elements $f\in m_D$, we have $\nu_D(f)\geq \nu_D(\pi'')$.  Let $s_{\epsilon}:=j\nu_D(\pi'')$.
 Note that $\pi''$ and $s_{\epsilon}$ both depend on the involution $\sigma$. 

\begin{lem}\label{150} Let $\epsilon \in \{1, -1\}$, $(D,\sigma)$ be a central $K-$division algebra as before, and assume that $s_{\epsilon}$ is odd. Then:

(a) There exists an $\epsilon-$symmetric uniformizer in $D$. 

(b) Any $\epsilon-$hermitian form over $(D,\sigma)$ is isometric to an $\epsilon-$hermitian form $h$ which can be decomposed as $h = h_0 \perp h_1$, such that any element in a diagonalization of $h_0$ is unit, and any element in a diagonalization of $h_1$ has value $\frac{1}{j}$.

\begin{proof}

(a): We have 
\begin{equation*} 
\begin{split}
\nu_D({\pi'}^{\frac{1-s_{\epsilon}}{2}} \pi'' \sigma({{\pi'}^{\frac{1-s_{\epsilon}}{2}}})) & = \frac{1-s_{\epsilon}}{2} \nu_D(\pi') + \nu_D(\pi'') + \frac{1-s_{\epsilon}}{2}\nu_D(\pi') \\
 & = \frac{1-s_{\epsilon}}{2j} + \frac{s_{\epsilon}}{j} + \frac{1-s_{\epsilon}}{2j} \\
 & = \frac{1}{j}.
\end{split}
\end{equation*}

We used the fact that the valuation $\nu_D$ is invariant under the involution $\sigma$ (see \cite[Page 14, Corollary 2.2]{4}). The element ${\pi'}^{\frac{1-s_{\epsilon}}{2}} \pi'' \sigma({{\pi'}^{\frac{1-s_{\epsilon}}{2}}})$ is an $\epsilon-$symmetric uniformizer. 

(b) Let $<u_1, \dots, u_r>$ be an $\epsilon-$hermitian form over $(D, \sigma)$. Replacing each $u_i$ by $\pi' u_i \sigma({\pi'})$ does not change the isometry class of the 1-dimensional $\epsilon-$hermitian form $<u_i>$, and the element $\pi' u_i \sigma({\pi'})$ is also $\epsilon-$symmetric, so 
$<u_i>\cong <\pi' u_i \sigma({\pi'})>$. On the other hand, $$\nu_D(\pi' u_i \sigma({\pi'}))=\nu_D (u_i)+\frac{2}{j}.$$

By replacing each $u_i$ by either $\pi' u_i \sigma({\pi'})$ or ${\pi'}^{-1} u_i \sigma({\pi'}^{-1})$ for suitable number of times, we can assume that each $u_i$ satisfies $\nu_D(u_i)\in \{0, \frac{1}{j}\}$, and the $K-$isometry class of $<u_1, \dots, u_r>$  does not change by these replacements. Hence we have the desired decomposition.  

\end{proof}
\end{lem}

\begin{thm} (Larmour)
If $h$ (from part (b) of Lemma \ref{150}) is $K-$anisotropic, then $h_0$ and $h_1$  are $K-$anisotropic, and this decomposition is unique up to isometry classes of $h_0$ and $h_1$. 
    \begin{proof}
        See \cite[Theorems 3.3, 3.4, and 3.6]{6}. 
    \end{proof}
\end{thm}

If $h\cong h_0$, we say that $h$ is unramified, and if $h\cong h_1$, we say that $h$ is ramified. We call $h=h_0 \perp h_1$ the \textbf{Larmour decomposition} of $h$, and we write $h=h_0 \perp h_1$ to represent this.\ \\

Let $\overline{D}:={\mathcal{O}_D}/{m_D}$ be the residue division ring. We will denote by $\overline{*}$ the image of $* \in \mathcal{O}_D$ in $\overline{D}$. The inclusion $\mathcal{O}\hookrightarrow \mathcal{O}_D$  induces an embedding $k={\mathcal{O}}/m \hookrightarrow \overline{D}= {\mathcal{O}_D}/{m_D}$. Elements in the image of $k$ under this embedding commute with all elements of $\overline{D}$, because $K$ is the center of $D$. Therefore, $\overline{D}$ is a $k-$algebra. So if $l$ is the center of $\overline{D}$, then $k$ embeds in $l$, and $\overline{D}$ will be an $l-$division algebra. The reduced norm map $Nrd: D \longrightarrow K$ is invariant under the involution $\sigma$ (see \cite[Page 14, Corollary 2.2]{4}). So $\sigma$ maps $\mathcal{O}_D$ to $\mathcal{O}_D$ and $m_D$ to $m_D$. The map

$$\overline{\sigma}: \overline{D} \longrightarrow \overline{D}$$
$$\  \   \  \   \    \ x + m_D \mapsto \sigma(x) + m_D$$
is an involution induced on $\overline{D}$.

 In order to have a similar isomorphism as Springer's theorem for $\epsilon-$hermitian forms over $(D,\tau)$, one needs to associate $h_0$ and $h_1$ some residue forms over $\overline{D}$. Let $h=h_0\perp h_1$ be  Larmour decomposition of $h$, $h_0=<u_1, \dots, u_p>$ and $h_1=<v_1, \dots, v_q>$. We denote by $\overline{h_0}$ the $\epsilon-$hermitian form $\overline{h_0}:=<\overline{u_1}, \overline{u_2}, \dots, \overline{u_p}>$ over $(\overline{D}, \overline{\sigma})$.

Let the integer $s_{\epsilon}$ be odd. Therefore, without loss of generality we may assume that $\pi'=\pi''$ (so $\pi'$ is $\epsilon-$symmetric), thanks to part (a) of the Lemma \ref{150}.

 Similarly, we want to associate a residue form to $h_1$. In order to do this, we need to consider twisted involutions. In general, for any element $\lambda\in D^*$, one can define the following twisted involution
$$\sigma_{\lambda}: D \longrightarrow D$$
$$\  \   \  \   \  \   \   \   \  \   \  \   \  \   \    \ x  \mapsto \lambda \sigma(x){\lambda}^{-1}.$$ Note that then $<v_1 {\pi'}^{-1}, \dots, v_q {\pi'}^{-1}>$ is a hermitian form over $(D,\sigma_{\pi'})$, which is unramified. We define the residue form of $h_1$ as $\overline{h_1}:=<\overline{v_1 {\pi'}^{-1}}, \dots, \overline{v_q {\pi'}^{-1}}>$. 

We denote the Witt group of hermitian forms over $(D,\sigma)$ by $W^1(D,\sigma)$, and the Witt group of skew-hermitian forms over $(D,\sigma)$ by $W^{-1}(D,\sigma)$.

\begin{thm}
    \label{21}

(Larmour) The map 

$$\partial= (\partial_0, \partial_1): W^{\epsilon}(D,\sigma) \longrightarrow  W^{\epsilon}(\overline{D},\overline{\sigma}) \oplus  W^1(\overline{D},\overline{\sigma_{\pi'}})$$
$$\   \  \ \ \   \  \    \   \    \   \  [h]   \     \ \mapsto  \  \   \   {([\overline{h_0}], [\overline{h_1}]})$$

is an isomorphism of groups.

\begin{proof}
See \cite[Theorem 3.7]{6}.
\end{proof}
\end{thm}

We will study Theorem \ref{21} in more details in Section \ref{52}, in the case of quaternion algebras. 

\begin{rema}\label{23}

If $s_{\epsilon}$ is even, then the same argument as in the proof of part (b) of Lemma \ref{150} shows that  every $\epsilon-$hermitian form $h$ over $(D,\tau)$ is unramified. Then Theorem \ref{21} gives the following isomorphism of groups:  

$$\partial=\partial_0: W^{\epsilon}(D,\sigma) \longrightarrow  W^{\epsilon}(\overline{D},\overline{\sigma})$$
$$\   \  \ \ \   \  \    \   \    \   \  [h]   \     \ \mapsto  \  \   \   [\overline{h}].$$
\end{rema}

\section{Larmour's decomposition for $\epsilon-$hermitian forms over quaternion algebras}\label{22}

We will continue to use the notations from the previous section, and from this point onward, we will assume that
$D$ is a quaternion division algebra over 
$K$. So $n=2$.  In this section we recall the different types of  ramification of $D$. Let $\tau$ be the canonical involution on $D$, and $x, y$ be two elements in $D$ such that  $x^2=a$, $y^2=b$, $D=(\frac{a,b}{K})$. Also put $z:=xy=-yx$. The set $\{1,x,y,z\}$ is a $K-$basis of $D$.

We have two cases of ramification for $D$: it is either unramified (i.e. $a$ and $b$ are units up to squares), or ramified (i.e. $a$ is a unit and $b$ is a uniformizer up to squares) (see \cite[Page 21, Example 1.17]{8}). The case where $D$ is generated by two uniformizers $a$ and $b$ (up to squares) is equivalent to the ramified case. This is because we have an isomorphism $(\frac{a,b}{K})\cong (\frac{a, -ab}{K})$ where $ab$ is a unit up to squares, and this isomorphism is due to the isometry of the norm forms associated to $(\frac{a,b}{K})$ and $(\frac{a, -ab}{K})$; they are both isometric to the form $<1,-a,-b,ab>$ (the isometry of the norm forms implies the isomorphism of the quaternion algebras, see \cite[Chapter 2, Theorem 2.5]{5}).\ \\

\textbf{Case A}: $D$ is unramified: $D=(\frac{a,b}{K})$ where both $a$ and  $b$ are units (so $j=1$). 

The assumption implies that $\pi$ is a uniformizer for $\nu_D$. We have $\Gamma_D=\mathbb{Z}$ and $\nu_D(x)=\nu_D(y)=\nu_D(z)=0$. Then in Larmour's decomposition $h= h_0 \perp h_1$, the form $h_0$ is unramified and in the diagonalization of $h_1$ each element has value $1$.

By \cite[Page 21, Example 1.17]{8}, $\overline{D}$ is isomorphic to $(\frac{\overline{a}, \overline{b}}{k})$, which is a division quaternion algebra (the norm form of $D$ is anisotropic and by Hensel's lemma, the norm form of $\overline{D}$ is anisotropic too, so $\overline{D}$ does not split). Now we compute $s_{\epsilon}$ and study the residue involution in the following subcases:\\

\textbf{Case A.1}: $\sigma$ is symplectic.

By \cite[Page 26, Proposition 2.21]{4}, $\sigma=\tau$.  Then the involution $\overline{\tau}$ is the following:

\begin{align*}
  \overline{\tau} : (\frac{\overline{a}, \overline{b}}{k}) & \rightarrow (\frac{\overline{a}, \overline{b}}{k}) \\
         e_0+ e_1 \overline{x} + e_2 \overline{y} + e_3 \overline{z} & \mapsto e_0 -e_1 \overline{x} - e_2 \overline{y} - e_3 \overline{z},
\end{align*}
for $e_0, e_1, e_2, e_3 \in k$.\\

The choice of $\pi''$ depends on $\epsilon$:

\textbf{Case A.1.1}:
If $\epsilon=1$, then we choose $\pi'=\pi''=\pi$. So $s_1=1$.

\textbf{Case A.1.2}:
If $\epsilon=-1$, then we choose $\pi'=\pi''=\pi x$. So $s_{-1}=1$. \\ \

\textbf{Case A.2}: $\sigma$ is  orthogonal.

Then $\sigma=\tau_{\zeta}$ where $\zeta\in Sym^{-1}(D,\tau)$ (by \cite[Page 26, Proposition 2.21]{4}), and the element $\zeta$ is uniquely determined by $\sigma$ up to a factor in $K^*$. Therefore, We may choose $\zeta$ from units of $\mathcal{O}_D$. Without loss of generality, we choose $\zeta=x$. The involution $\overline{\tau_x}$ is the following map:

\begin{align*}
  \overline{\tau_x} : (\frac{\overline{a}, \overline{b}}{k}) & \rightarrow (\frac{\overline{a}, \overline{b}}{k}) \\
         e_0+ e_1 \overline{x} + e_2 \overline{y} + e_3 \overline{z} & \mapsto e_0 -e_1 \overline{x} + e_2 \overline{y} + e_3 \overline{z},
\end{align*}
for $e_0, e_1, e_2, e_3 \in k$.\\

\textbf{Case A.2.1}:
If $\epsilon=1$, then we choose $\pi'=\pi''=\pi$. So $s_1=1$.

\textbf{Case A.2.2}:
If $\epsilon=-1$, then we choose $\pi'=\pi''=\pi x$. So $s_{-1}=1$. \\ \

\textbf{Case B}: $D$ is ramified: $D=( \frac{a, \pi}{K})$ where $a$ is a unit (so $j=2$). In fact if $D=(\frac{a,b}{K})$ with $\nu_K(a)=0$ and $\nu_K(b)=1$, then we can choose $b$ to be the uniformizer of $K$, and hence without loss of generality assume that $b=\pi$.

In this case $y=\sqrt{\pi}$  is a uniformizer for $\nu_D$. We have $\Gamma_D= \{\frac{l}{2} \ | \ l \in \mathbb{Z}\}$, $\nu_D(x)=0$, and $\nu_D(y)=\nu_D(z)=\frac{1}{2}$. In Larmour's decomposition $h= h_0 \perp h_1$, the form $h_0$ is unramified and in the diagonalization of $h_1$ each element has value $\frac{1}{2}$. 
By \cite[Page 21, Example 1.17]{8}, $\overline{D}$ is isomorphic to $k(\overline{x})$, which is a quadratic field extension of $k$. Note that in general, $k(\overline{x})$ is a quadratic \'{e}tale extension of $k$, but in our case we know that $\overline{x}$ is not a square in $k$ since otherwise by Hensel's lemma $x$ will be a square in $K$ which contradicts the fact that $D$ is a division quaternion algebra. Therefore, $k(\overline{x})$ is a quadratic field extension of $k$. Let $\iota$ denote the nontrivial $k-$automorphism of $k(\overline{x})$.

\begin{align*}
  \iota:  k(\overline{x}) & \rightarrow k(\overline{x})  \\
         e_0 + e_1 \overline{x} & \mapsto e_0 - e_1 \overline{x},
\end{align*} for $e_0 , e_1 \in k$.

Now we compute $s_{\epsilon}$ and study the residue involution in the following subcases:\\

\textbf{Case B.1}: $\sigma$ is symplectic (hence $\sigma=\tau$). 

We have $\overline{\sigma}=\iota$ which is an involution of the second kind. Now we choose $\pi''$ and compute $s_{\epsilon}$:

\textbf{Case B.1.1}: If $\epsilon=1$, then we choose $\pi'=y$ and  $\pi''=\pi$. So $s_1=2$.

\textbf{Case B.1.2}: If $\epsilon=-1$, then we choose $\pi'=\pi''=y$. So $s_{-1}=1$. \\ \

\textbf{Case B.2}: $\sigma$ is orthogonal. 

Similar to case A.2, $\sigma=\tau_{\zeta}$ where $\zeta\in Sym^{-1}(D,\tau)$, and the element $\zeta$ is uniquely determined by $\sigma$ up to a factor in $K^*$. So value $\nu_D(\zeta)$ can assumed to be $0$ or $\frac{1}{2}$. 

         \textbf{Case B.2.1}: $\nu_D(\zeta)=0$. Without loss of generality, we may assume $x=\zeta$. Choose a vector $w$ from the orthogonal complement of the vector subspace spanned by $1$ and $x$, i.e. from $<1,x>^{\perp}$. By a similar argument as in the proof of the first part of Lemma \ref{150}, we can assume that $\nu_D(w)\in \{0, \frac{1}{2}\}$. However, $\nu_D(w)$ cannot be $0$, because in that case $D$ will be unramified, which contradicts our assumption. Hence we have $\nu_D(w)=\frac{1}{2}$. Without loss of generality assume that $y=w$, and then $z=\zeta w$. So $Sym^1(D,\tau_{\zeta})=<1,y,z>$ and $Sym^{-1}(D,\tau_{\zeta})=<x>$. 
        
The involution $\overline{\tau_{\zeta}}$ will be the map $\iota$, similar to case B.1 (which is an involution of the second kind):

\begin{align*}
 \iota= \overline{\tau_{\zeta}} :  k(\overline{x}) & \rightarrow k(\overline{x})  \\
         e_0 + e_1 \overline{x} & \mapsto e_0 - e_1 \overline{x},
\end{align*} for $e_0 , e_1 \in k$.

Case B.2.1 splits into two subcases:

\textbf{Case B.2.1.1}: If $\epsilon=1$, then we choose $\pi'=\pi''=y$. So $s_1=1$.

\textbf{Case B.2.1.2}: If $\epsilon=-1$, then we choose $\pi'=y$ and $\pi''=\pi x$. So $s_{-1}=2$. \\ \

  \textbf{Case B.2.2}: $\nu_D(\zeta)=\frac{1}{2}$.
  
  Without loss of generality, we may assume $y=\zeta$. Choose a vector $w'$ from the orthogonal complement of the vector subspace spanned by $1$ and $y$, i.e. from $<1,y>^{\perp}$. By a similar argument as in the proof of the first part of Lemma \ref{150}, we can assume that $\nu_D(w')\in \{0, \frac{1}{2}\}$. 
   \begin{itemize}
       \item \rom{1}: If $\nu_D(w')=0$, then without loss of generality assume that $x=w'$ and $z=w'\zeta$.
       \item  \rom{2}: If $\nu_D(w')=\frac{1}{2}$, then without loss of generality assume that $z=w'$ and $x=w' {\zeta}^{-1}$.
   \end{itemize}
    So in any case, $x$ is a unit, and $z$ is a uniformizer. Hence the two possibilities \rom{1} and \rom{2} above are equivalent, and without loss of generality, in case B.2.2 we assume that $x=w'$ and $z=w'\zeta$ (the case \rom{1} above). We have $Sym^1(D,\tau_{\zeta})=<1,x,z>$ and $Sym^{-1}(D,\tau_{\zeta})=<y>$.
    
   Since $K(x)\subseteq Sym^1(D, \tau_{\zeta})$, the involution $\overline{\tau_{\zeta}}$ will exactly be the identity map of the field $k(\overline{x})$.

\begin{align*}
  \overline{\tau_{\zeta}} :  k(\overline{x}) & \rightarrow k(\overline{x})  \\
         e_0 + e_1 \overline{x} & \mapsto e_0 + e_1 \overline{x},
\end{align*} for $e_0 , e_1 \in k$.

We have two subcases:

\textbf{Case B.2.2.1}: If $\epsilon=1$, then we choose $\pi'=\pi''=z$. So $s_1=1$.

\textbf{Case B.2.2.2}: If $\epsilon=-1$, then we choose $\pi'=\pi''=y$. So $s_{-1}=1$. \ \\

We summarize all cases in the table below 
(QDA in column 5 referes to the quaternion division algebra $(\frac{\overline{a}, \overline{b}}{k})$):

\begin{equation}\label{45}
    \begin{tabular}{||c  c c c c c c c c ||} 
 \hline
 Case & $j$ & $\pi'$ & $\sigma$ & $\overline{D}$ & $\overline{\sigma}$ & $\epsilon$ & $\pi''$ & $s_{\epsilon}$    \\ [0.8ex] 
 \hline\hline
 A.1.1 & $1$ & $\pi$ & $\tau$   & QDA  & $\overline{\tau}$     & $1$ & $\pi$ & $1$  \\ 
 \hline
  A.1.2 & $1$ & $\pi x$ & $\tau$  & QDA & $\overline{\tau}$ & $-1$ & $\pi x$ & $1$  \\
 \hline
  A.2.1 & $1$ & $\pi$ & $\tau_x$  & $k(\overline{x})$ & $\iota$ & $1$ & $\pi$ & $1$  \\
 \hline
  A.2.2 & $1$ & $\pi x$ & $\tau_x$ & $k(\overline{x})$ & $\iota$  & $-1$ & $\pi x$ & $1$  \\
 \hline
  B.1.1 & $2$ & $y$ & $\tau$ & 
$k(\overline{x})$ & $\iota$ & $1$ & $\pi$ & $2$ \\
 \hline
  B.1.2 & $2$ & $y$ & $\tau$ & $k(\overline{x})$ & $\iota$  & $-1$ &  $y$ & $1$  \\
 \hline
  B.2.1.1 & $2$ & $y$ & $\tau_x$ & $k(\overline{x})$ & $\iota$  & $1$ & $y$ & $1$  \\
 \hline
   B.2.1.2 & $2$ & $y$ & $\tau_x$ & $k(\overline{x})$ & $\iota$ & $-1$ & $\pi x$ & $2$    \\ \hline 
   B.2.2.1 & $2$ & $z$ & $\tau_y$ & $k(\overline{x})$ & id & $1$ & $z$ & $1$  \\
 \hline

  B.2.2.2 & $2$  & $y$ & $\tau_y$ & $k(\overline{x})$ & id & $-1$ & $y$ & $1$   \\ \hline

 \hline
\end{tabular}
\end{equation}

\begin{rema}
    
 Larmour's decomposition $h= h_0 \perp h_1$ for an $\epsilon-$hermitian form $h$ is trivial in cases 2.1.1 and 2.2.1.2, because $h$ will be automatically unramified (equivalently, $s_{\epsilon}=2$).  
\end{rema}

\section{Diagonal entries in Larmour decomposition}\label{20}

Recall that $h$ is a $K-$anisotropic $\epsilon-$hermitian form over $(D,\sigma)$ where $\sigma=\tau$, or $\sigma=\tau_{\zeta}$ for $\zeta\in \{x,y\}$. Let $h=h_0 \perp h_1$ be Larmour decomposition of $h$. In this section, we describe how the diagonal entries in $h_0$ and $h_1$ look like.

  Any element $u$ in the diagonalization of $h$ belongs to $Sym^{\epsilon}(D,\sigma)$. Since $D$ is spanned by $1,x,y,z$, there exist elements $\delta, \alpha, \beta, \gamma \in K$ such that $u=\delta + \alpha x + \beta y + \gamma z$. A priori, it seems that the  elements $\delta, \alpha, \beta, \gamma$ need not necessarily be in $\mathcal{O}_K$. We will show that for any $h$, up to isometry in $K$, there exist a diagonalization such that the elements $\delta, \alpha, \beta, \gamma$ belong to $\mathcal{O}_K$. We will also prove that $u$ belongs to a specific subspace of $D$ in each case, depending on whether $u$ is a diagonal entry of $h_0$ or $h_1$.

\begin{lem}\label{47}
  The ring ${\mathcal{O}_D}$ is contained in the $\mathcal{O}_K$- submodule of $D$ generated by elements $1,x,y,z$.  
   \begin{proof}
   Let $u\in{\mathcal{O}_D}$.  
       The division algebra $D$ is generated (as a $K-$vector space) by elements $1,x,y,z$, so there exist elements $\delta, \alpha, \beta, \gamma \in K$ such that $u=\delta + \alpha x + \beta y + \gamma z$. We need to show that $\delta, \alpha, \beta, \gamma\in \mathcal{O}_K$. Without loss of generality, assume that $u\ne0$. Let $m=min(\nu_K(\delta), \nu_K(\alpha), \nu_K(\beta), \nu_K(\gamma))$. We want to show that $m$ is nonnegative.  Assume the contrary, that $m<0$. 

       Let 
      \begin{equation*}
\delta'=
    \begin{cases}
        \pi^{-m}\delta &\delta\ne 0\\
        0 & \delta=0
    \end{cases}
\end{equation*} and similarly define  $\alpha', \beta', \gamma'$. Note that elements $\delta, \alpha, \beta, \gamma$ cannot be simultaneously zero. We have

\begin{equation*} 
\begin{split}
0 & = \nu_D(u) \\
 & = \frac{1}{2}\nu_K(Nrd(u))\\
 & = \frac{1}{2} \nu_K(u \tau(u)) \\
 & = \frac{1}{2}\nu_K(\delta^2 -  \alpha^2 a - \beta^2 b + \gamma^2 ab) \\
 & = \frac{1}{2}(-2m+\nu_K(\delta'^2 -  \alpha'^2 a - \beta'^2 b + \gamma'^2 ab)) \\
 & = -m + \frac{1}{2} \nu_K(\delta'^2 -  \alpha'^2 a - \beta'^2 b + \gamma'^2 ab). 
\end{split}
\end{equation*}

Therefore, $\delta'^2 -  \alpha'^2 a - \beta'^2 b + \gamma'^2 ab \in m_K$. Consider the norm form $q=<1, -a, -b, ab>$ over $K$.

If $D$ is unramified, then consider the quadratic form $\overline{q}=<\overline{1}, -\overline{a}, -\overline{b}, \overline{ab}>$ over $k$. Then $\overline{q}$ will be $k-$ isotropic because  $\overline{q}(\overline{\delta'}, \overline{\alpha'}, \overline{\beta'}, \overline{\gamma'})=0$. By Springer's decomposition theorem for quadratic forms, $q$ will be $K-$isotropic. Hence $D$ will be split over $K$, which is a contradiction. Therefore, 
$\delta, \alpha, \beta, \gamma\in \mathcal{O}_K$.

If $D$ is ramified, then $b, ab \in m_K$, hence $-\beta'^2 b  +  \gamma'^2 ab\in m_K$. Consider the quadratic form $\overline{q}=<\overline{1}, -\overline{a}>$ over $k$. Then $\overline{q}$ will be $k-$isotropic because $\overline{q}(\overline{\delta'}, \overline{\alpha'})=0$. By Springer's decomposition theorem for quadratic forms, $q$ will be $K-$isotropic. Hence $D$ will be split over $K$, which is a contradiction. Therefore 
$\delta, \alpha, \beta, \gamma\in \mathcal{O}_K$.
   \end{proof}
\end{lem}

Let $u$ be an element in the diagonalization of the $\epsilon-$hermitian form $h$ over $(D,\sigma)$. By Lemma \ref{47}, $u$ belongs to the $\mathcal{O}_K-$submodule of $D$ generated by the elements $\{1,x,y,z\}\cap Sym^{\epsilon}(D,\sigma)$. When $D$ is ramified, then we will have  a stronger description of the diagonal entries: we will show that $u$ belongs to the $\mathcal{O}_K-$submodule of $D$ generated by the elements $\{1,x\}\cap Sym^{\epsilon}(D,\sigma)$ if $u$ is from the unramified part (i.e. $h_0$), and $u$ belongs to the $\mathcal{O}_K-$submodule of $D$ generated by the elements $\{y,z\}\cap Sym^{\epsilon}(D,\sigma)$ if $u$ is from the ramified part (i.e. $h_1$).

 The following lemma will be helpful in our computations.

\begin{lem}
    \label{64}
Let $v_0,v_1\in D$ be symmetric or skew-symmetric units such that $\overline{v_0}=\overline{\sigma}(\theta) \overline{v_1} \theta$ for some $\theta \in \overline{D}$. Then there is $t\in D$ with $\overline{t}=\theta$ and $v_0= \sigma(t) v_1 t$.
    \begin{proof}
    See \cite[Lemma 2.2]{6}.
    \end{proof}
\end{lem}

\begin{lem}\label{40}
    Let $D$ be ramified (equivalently, $j=2$), and $h=<\delta + \alpha x + \beta y +\gamma z>$ be a one-dimensional unramified $\epsilon-$hermitian form over $(D,\sigma)$, where $\delta,  \alpha,  \beta, \gamma \in \mathcal{O}_K$. Then $h\cong <\delta + \alpha x>$.
    \begin{proof}
       Since $y,z \in m_D$, we have $\overline{\delta + \alpha x + \beta y +\gamma z}= \overline{\delta + \alpha x}$. Applying Lemma \ref{64} with $v_0:=\delta + \alpha x + \beta y +\gamma z$, $v_1= \delta + \alpha x$, and $\theta=\overline{1}$, we get an element $t\in D$ such that $\delta + \alpha x + \beta y +\gamma z = \sigma(t) (\delta + \alpha x) t$, which implies that $h\cong <\delta + \alpha x>$.
    \end{proof}
\end{lem}

\begin{rema}\label{41}
    Let $D$ be ramified (equivalently, $j=2$), and $h=<\delta + \alpha x + \beta y +\gamma z>$ be a one-dimensional unramified $\epsilon-$hermitian form over $(D,\sigma)$, where $\delta,  \alpha,  \beta, \gamma \in \mathcal{O}_K$. Then Lemma \ref{40} implies that:  

(1) If $\sigma=\tau$ and $\epsilon=1$ (case B.1.1), then $h\cong <\delta>$ and $\delta \in {\mathcal{O}_K}^*$.

(2) If $\sigma=\tau$ and $\epsilon=-1$ (case B.1.2), then $h\cong <\alpha x>$ and $\alpha \in {\mathcal{O}_K}^*$.

(3) If $\sigma=\tau_x$ and $\epsilon=1$ (case B.2.1.1), then $h\cong <\delta>$ and $\delta \in {\mathcal{O}_K}^*$.

(4) If $\sigma=\tau_x$ and $\epsilon=-1$ (case B.2.1.2), then $h\cong <\alpha x>$ and $\alpha \in {\mathcal{O}_K}^*$.

(5) If $\sigma=\tau_y$, then $\epsilon$ has to be $1$ (case B.2.2.1), $h\cong <\delta + \alpha x>$ and $\delta + \alpha x$ is a unit in the field $K(x)$.
\end{rema}

Now we investigate the case where $h$ is ramified.

\begin{lem}\label{42}
Let $D$ be ramified (equivalently, $j=2$), $h=<\delta + \alpha x + \beta y +\gamma z>$ be a one-dimensional ramified $\epsilon-$hermitian form over $(D,\sigma)$, where $\delta,  \alpha,  \beta, \gamma \in \mathcal{O}_K$. Then $h\cong <\beta y + \gamma z>$.

\begin{proof} Let $u=\delta + \alpha x + \beta y +\gamma z$. The extension $K(x)/K$ is unramified, so the value group of $K(x)$ is same as the value group of $K$, which is $\mathbb{Z}$. Hence $\nu_D(\delta+\alpha x)\in \mathbb{Z}$. Since $\nu_D(u)=\nu_D(\delta + \alpha x +  \beta y + \gamma z)=\frac{1}{2}$, we have 

$$\nu_D(\delta +\alpha x ) > \frac{1}{2}   \   \ \ \  \ \  \    \textit{and}   \   \ \ \  \ \  \ \nu_D(\beta y + \gamma z )=\frac{1}{2}, $$
which forces $\nu_D(\delta + \alpha x)$ to be at least $1$. The element $(\beta y + \gamma z)y^{-1}= \beta + \gamma x$ will then be a unit. 

In cases B.1.2, B.2.1.1, and B.2.2.2, $y\in Sym^1(D,\sigma)$, hence $uy^{-1}\in Sym^1(D,\sigma_y)$. Also in case B.2.2.1,  $y\in Sym^{-1}(D,\sigma)$, hence $uy^{-1}\in Sym^{-1}(D,\sigma_y)$. Furthermore, $uy^{-1}$ is a unit in all cases. The fact that $\nu_D(\delta + \alpha x)\geq 1$ implies that $\nu_D((\delta + \alpha x)y^{-1})\geq \frac{1}{2}$, hence $(\delta + \alpha x)y^{-1} = \pi^{-1}\delta y - \alpha z \in m_D$. Therefore, $\overline{uy^{-1}}=\overline{\beta + \gamma x + \pi^{-1}\delta y - \alpha z}=\overline{ \beta + \gamma x}$, so by applying Lemma \ref{64} to $\sigma=\sigma_y$, $v_0= \beta + \gamma x + \pi^{-1}\delta y - \alpha z$, $v_1=\beta + \gamma x$ and $\theta=\overline{1}$, we get an element $t \in D$ such that $uy^{-1}=\beta + \gamma x + \pi^{-1}\delta y - \alpha z = \sigma_y(t) (\beta + \gamma x) t$. Hence $<uy^{-1}>=<\beta + \gamma x + \pi^{-1}\delta y - \alpha z> \cong <\beta + \gamma x>$ as $\epsilon-$hermitian forms over $(D,\sigma_y)$, which implies that $h=<u>=<\delta + \alpha x +  \beta y + \gamma z> \cong <\beta y + \gamma z>$
as $\epsilon-$hermitian forms over $(D,\sigma)$.
\end{proof}
\end{lem}

\begin{rema}\label{43}

Let $D$ be ramified (equivalently, $j=2$), $h=<\delta + \alpha x + \beta y +\gamma z>$ be a one-dimensional ramified $\epsilon-$hermitian form over $(D,\sigma)$, where $\delta,  \alpha,  \beta, \gamma \in \mathcal{O}_K$. Then  Lemma \ref{42} implies that:

(1) If $\sigma=\tau$ then $\epsilon$ has to be $-1$ (case B.1.2), $h\cong <\beta y + \gamma z>$ and $\beta + \gamma x$ is a unit in the field $K(x)$.

(2) If $\sigma=\tau_x$ then $\epsilon$ has to be $1$ (case B.2.1.1), $h\cong <\beta y + \gamma z>$ and $\beta + \gamma x$ is a unit in the field $K(x)$.

(3) If $\sigma=\tau_y$ and $\epsilon=1$ (case B.2.2.1), then $h\cong <\gamma z>$ and $\gamma \in {\mathcal{O}_K}^*$.

(4) If $\sigma=\tau_y$ and $\epsilon=-1$ (case B.2.2.2), then  $h\cong <\beta y>$ and $\beta\in {\mathcal{O}_K}^*$.
\end{rema}

The results in this section are summarized in the following table:

\begin{equation}\label{49}
\begin{tabular}{||c  c c c c c c c c c ||} 
 \hline
 Case & $j$ & $\pi'$ & $\sigma$ & $\epsilon$ & $\pi''$ & $s_{\epsilon}$ & $Sym^{\epsilon}(D,\sigma)$ & $h_0$ & $h_1$    \\ [0.8ex] 
 \hline\hline
 A.1.1 & $1$ & $\pi$ & $\tau$   &  $1$ & $\pi$ & $1$ & $<1>$ & $\delta$ & $\delta$   \\ 
 \hline
  A.1.2 & $1$ & $\pi x$ & $\tau$  &  $-1$ & $\pi x$ & $1$ & $<x,y,z>$ & $\alpha x + \beta y + \gamma z$ & $\alpha x + \beta y + \gamma z$  \\
 \hline
  A.2.1 & $1$ & $\pi$ & $\tau_x$  &  $1$ & $\pi$ & $1$ & $<1,y,z>$ & $\delta + \beta y + \gamma z$ & $\delta + \beta y + \gamma z$ \\
 \hline
  A.2.2 & $1$ & $\pi x$ & $\tau_x$ &  $-1$ & $\pi x$ & $1$ & $<x>$ & $\alpha x$ & $\alpha x$ \\
 \hline
  B.1.1 & $2$ & $y$ & $\tau$ & 
 $1$ & $\pi$ & $2$ & $<1>$ & $\delta$ & -  \\
 \hline
  B.1.2 & $2$ & $y$ & $\tau$ &  $-1$ &  $y$ & $1$ & $<x,y,z>$ & $\alpha x$ & $\beta y + \gamma z$  \\
 \hline
  B.2.1.1 & $2$ & $y$ & $\tau_x$ &  $1$ & $y$ & $1$ & $<1,y,z>$ & $\delta$ & $\beta y + \gamma z$ \\
 \hline
   B.2.1.2 & $2$ & $y$ & $\tau_x$ &  $-1$ & $\pi x$ & $2$ & $<x>$ & $\alpha x$ & -   \\ \hline 
   B.2.2.1 & $2$ & $z$ & $\tau_y$ &  $1$ & $z$ & $1$ & $<1,x,z>$ & $\delta + \alpha x$ & $\gamma z$ \\
 \hline
  B.2.2.2 & $2$  & $y$ & $\tau_y$ &  $-1$ & $y$ & $1$ & $<y>$ & - & $\beta y$  \\ \hline 
\end{tabular}
\end{equation}

\section{Explicit formulas for the isomorphism of Witt groups}\label{52}

Recall the isomorphism (see Theorem \ref{21}):

$$\partial= (\partial_0, \partial_1): W^{\epsilon}(D,\sigma) \longrightarrow  W^{\epsilon}(\overline{D},\overline{\sigma}) \oplus  W^1(\overline{D},\overline{\sigma_{\pi'}})$$
$$\   \  \ \ \   \  \    \   \    \   \  [h]   \     \ \mapsto  \  \   \   {([\overline{h_0}], [\overline{h_1}]}).$$

Now we consider the cases introduced in Section \ref{22} again, and for each case we describe the above isomorphism  explicitly and compute the residue forms $\overline{h_0}$ and $\overline{h}_1$ for one-dimensional forms $h_0$ and $h_1$ (in what follows, by "an $\epsilon-$hermitian form $h$" we actually mean "the $K-$isometry class of an $\epsilon-$hermitian form $h$"). Since the Witt groups are generated by one dimensional forms, 
in order to have an explicit formula for the maps $\partial_0$ and $\partial_1$, it is enough to compute the image for one-dimensional forms $h_0$ and $h_1$. If $h_i=<u>$ for $i \in \{0,1\}$, then by Lemma \ref{47}, we may assume that $u=\delta + \alpha x + \beta y + \gamma z$ for some $\delta, \alpha, \beta, \gamma \in \mathcal{O}_K$, such that

      \begin{equation*}
\nu_D(u)=
    \begin{cases}
        0 & i=0\\
        1 & i=1 \ \textit{and} \ $D$  \ \textit{is unramified}\\
        \frac{1}{2} & i=1 \ \textit{and} \  $D$ \  \textit{is ramified}\\
    \end{cases}
\end{equation*} \\

\textbf{Case A}: $D$ is unramified.

Recall that in this case $D=(\frac{a,b}{K})$ where both $a$ and  $b$ are units, $j=1$, and $\overline{D}$ is isomorphic to $(\frac{\overline{a}, \overline{b}}{k})$, which is a division quaternion algebra over $k$.

\textbf{Case A.1}: $\sigma=\tau$ is symplectic.

Hence $\overline{\sigma}=\overline{\tau}$ is the symplectic involution on $\overline{D}$. 

\textbf{Case A.1.1}: $\epsilon=1$. 

We have $\pi'=\pi\in K$. Hence $\sigma_{\pi'}=\tau_{\pi'}=\tau$, which gives the following isomorphism:\\ 

$$W^{1}((\frac{a,b}{K}),\tau) \cong W^{1}((\frac{\overline{a},\overline{b}}{k}),\overline{\tau}) \oplus W^{1}((\frac{\overline{a},\overline{b}}{k}),\overline{\tau}).$$

Note that $Sym^1((\frac{a,b}{K}),\tau)=K$ and $Sym^1((\frac{\overline{a},\overline{b}}{k}),\overline{\tau})=k$. The image of any one dimensional unramified hermitian form $<\delta>$ (where $\delta\in {\mathcal{O}_K}^*$) under the 
 map $\partial_0$ is $<\overline{\delta}>$. The image of any one dimensional ramified hermitian form $<\delta>$ (where $\nu_K(\delta)=1$) under the 
 map $\partial_1$ is $<\overline{\pi^{-1}\delta}>$.

\textbf{Case A.1.2}: $\epsilon=-1$.

We have $\pi'=\pi x$. Hence $\sigma_{\pi'}=\tau_{\pi'}=\tau_{x}$ is an orthogonal involution. We have the following isomorphism:\\ 

$$W^{-1}((\frac{a,b}{K}),\tau) \cong W^{-1}((\frac{\overline{a},\overline{b}}{k}),\overline{\tau}) \oplus W^{1}((\frac{\overline{a},\overline{b}}{k}),\overline{\tau_x}).$$

Note that $Sym^{-1}((\frac{a,b}{K}),\tau)=<x,y,z>$,  $Sym^{-1}((\frac{\overline{a},\overline{b}}{k}),\overline{\tau})=<\overline{x}, \overline{y}, \overline{z}>$, and $Sym^1((\frac{\overline{a},\overline{b}}{k}),\overline{\tau_x})=<\overline{1}, \overline{y}, \overline{z}>$. The image of any one dimensional unramified skew-hermitian form $<\alpha x + \beta y + \gamma z>$ (where $\alpha, \beta, \gamma\in \mathcal{O}_K$ and $\alpha x + \beta y + \gamma z\in {\mathcal{O}_D}^*$) under the 
 map $\partial_0$ is $<\overline{\alpha x + \beta y + \gamma z}>$. 
The image of any one dimensional ramified skew-hermitian form $<\alpha x + \beta y + \gamma z>$ (where $\alpha, \beta, \gamma \in m_K$, and $\nu_D(\alpha x + \beta y + \gamma z)=1$) under the 
 map $\partial_1$ is $<\overline{\pi^{-1}(\alpha  -\gamma y + \beta z)}>$. 

\textbf{Case A.2}: $\sigma=\tau_x$ is orthogonal.

Hence $\overline{\sigma}=\overline{\tau_x}$ is an orthogonal involution on $\overline{D}$. 

\textbf{Case A.2.1}: $\epsilon=1$. 

We have $\pi'=\pi\in K$. Hence $\sigma_{\pi'}={\tau_x}_{\pi'}={\tau_x}_{\pi}=\tau_x$, which gives the following isomorphism:\\ 

$$W^{1}((\frac{a,b}{K}),\tau_x) \cong W^{1}((\frac{\overline{a},\overline{b}}{k}),\overline{\tau_x}) \oplus W^{1}((\frac{\overline{a},\overline{b}}{k}),\overline{\tau_x}).$$

Note that $Sym^1((\frac{a,b}{K}),\tau_x)=<1,y,z>$ and $Sym^1((\frac{\overline{a},\overline{b}}{k}),\overline{\tau_x})=<\overline{1}, \overline{y}, \overline{z}>$. The image of any one dimensional unramified hermitian form $<\delta + \beta y + \gamma z>$ (where $\delta,\beta,\gamma\in {\mathcal{O}_K}$, and $\delta + \beta y + \gamma z \in {\mathcal{O}_D}^*$) under the 
 map $\partial_0$ is $<\overline{\delta + \beta y + \gamma z}>$. The image of any one dimensional ramified hermitian form $<\delta + \beta y + \gamma z>$ (where $\delta,\beta,\gamma\in m_K$ and $\nu_D(\delta + \beta y + \gamma z)=1$) under the 
 map $\partial_1$ is $<\overline{\pi^{-1}(\delta + \beta y + \gamma z)}>$.

\textbf{Case A.2.2}: $\epsilon=-1$.

We have $\pi'=\pi x$. Hence $\sigma_{\pi'}={\tau_x}_{\pi'}=\tau_{\pi x^2}=\tau_{a\pi}=\tau$ is the symplectic  involution. We have the following isomorphism:\\ 

$$W^{-1}((\frac{a,b}{K}),\tau_x) \cong W^{-1}((\frac{\overline{a},\overline{b}}{k}),\overline{\tau_x}) \oplus W^{1}((\frac{\overline{a},\overline{b}}{k}),\overline{\tau}).$$

Note that $Sym^{-1}((\frac{a,b}{K}),\tau_x)=<x>$,  $Sym^{-1}((\frac{\overline{a},\overline{b}}{k}),\overline{\tau_x})=<\overline{x}>$, and $Sym^1((\frac{\overline{a},\overline{b}}{k}),\overline{\tau})=k$. The image of any one dimensional unramified skew-hermitian form $<\alpha x>$ (where $\alpha \in {\mathcal{O}_K}^*$) under the 
 map $\partial_0$ is $<\overline{\alpha} 
 \ \overline{x}>$. 
The image of any one dimensional ramified skew-hermitian form $<\alpha x>$ (where $\nu_K(\alpha)=1$) under the 
 map $\partial_1$ is $<\overline{\pi^{-1} \alpha}>$.

\textbf{Case B}: $D$ is ramified.

Recall that in this case $D=(\frac{a, \pi}{K})$ where $a$ is a unit, $j=2$, and  $\overline{D}$ is isomorphic to $k(\overline{x})$ (which is a quadratic field extension of $k$), and the involution $\overline{\tau}$ is the nontrivial $k-$automorphism of $k(\overline{x})$, which we denote by  $\iota$. In case B, the Witt group of quadratic forms over the fields $K(x)$ and $k(\overline{x})$ ($W(K(x))$ and $W(k(\overline{x}))$) will play an important role.
 
\textbf{Case B.1}:  $\sigma=\tau$ is symplectic. 

Hence $\overline{\sigma}=\overline{\tau}=\iota$. 

\textbf{Case B.1.1}: $\epsilon=1$. 

We have $s_{\epsilon}=2$ and by Remark \ref{23} we have $$W^{1}((\frac{a,\pi}{K}),\tau) \cong W^1(k(\overline{x}),\iota).$$

Note that $Sym^{1}((\frac{a,\pi}{K}),\tau)=K$ and $Sym^{1}(k(\overline{x}), \iota)=k$. The image of any one dimensional hermitian form $<\delta>$ which is automatically unramified ($\alpha \in {\mathcal{O}_K}^*$) under the 
 map $\partial_0$ is $<\overline{\delta}>$. 
 
\textbf{Case B.1.2}: $\epsilon=-1$.

We have $\pi'=y$. Hence $\sigma_{\pi'}=\sigma_y=\tau_y$ is an orthogonal involution. Since $K(x)\subseteq Sym^{1}(D,\tau_y)=<1,x,z>$, the restriction of the involution $\tau_y$ on $K(x)$ is trivial. Therefore $\overline{\sigma_{\pi'}}=id$ so $W^1(k(\overline{x}),\tau_y)\cong W(k(x))$. Hence 

$$W^{-1}((\frac{a,\pi}{K}),\tau) \cong  W^{-1}(k(\overline{x}),\iota) \oplus  W(k(\overline{x})).$$

Note that $Sym^{-1}((\frac{a,\pi}{K}),\tau)=<x,y,z>$ and $Sym^{-1}(k(\overline{x}),\iota)=<\overline{x}>$. By Remark \ref{41} part (2), any one dimensional unramified skew-hermitian form can be diagonalized as $<\alpha x>$ for $\alpha \in {\mathcal{O}_K}^*$. The image of this form under $\partial_0$ is $<\overline{\alpha} \ \overline{x}>$. Also by Remark \ref{43} part (1), any one dimensional ramified skew-hermitian form  can be diagonalized as $<\beta y + \gamma z>$ where $\beta + \gamma x$ is a unit in the field $K(x)$. The image of this form under $\partial_1$ is the quadratic form $<\overline{\beta}  + \overline{\gamma} \ \overline{x}>$.

\textbf{Case B.2}: $\sigma=\tau_{\zeta}$ is orthogonal, where $\zeta \in \{x,y\}$.

\textbf{Case B.2.1}: $\zeta=x$. 

\textbf{Case B.2.1.1}: $\epsilon=1$.

We have $\pi'=y$. Hence $\sigma_{\pi'}=\sigma_y={\tau_{x}}_y=\tau_z$ is an orthogonal involution. Since $K(x)\subseteq Sym^1(D,\tau_z)=<1,x,y>$, the restriction of the involution $\tau_z$ on $K(x)$ is trivial. Therefore, $\overline{\sigma_{\pi'}}=id$ so $W^1(k(\overline{x}),\overline{\tau_z})\cong W(k(x))$. Hence

$$W^{1}((\frac{a,\pi}{K}),\tau_x) \cong  W^{1}(k(\overline{x}),\iota) \oplus  W(k(\overline{x})).$$

Note that $Sym^{1}((\frac{a,\pi}{K}),\tau_x)=<1,y,z>$ and $Sym^{1}(k(\overline{x}),\iota)=k$. By Remark \ref{41} part (3), any one dimensional unramified hermitian form can be diagonalized as $<\delta>$ for $\delta \in {\mathcal{O}_K}^*$. The image of this form under $\partial_0$ is $<\overline{\delta}>$. Also by Remark \ref{43} part (2), any one dimensional ramified hermitian form can be diagonalized as $<\beta y + \gamma z>$ where $\beta + \gamma x$ is a unit in the field $K(x)$. The image of this form under $\partial_1$ is the quadratic form  $<\overline{\beta}  + \overline{\gamma} \ \overline{x}>$.

\textbf{Case B.2.1.2}: $\epsilon=-1$.

We have $s_{\epsilon}=2$ and by Remark \ref{23} we have $$W^{-1}((\frac{a,\pi}{K}),\tau_x) \cong W^{-1}(k(\overline{x}),\iota).$$

Note that $Sym^{-1}((\frac{a,\pi}{K}),\tau_x)=<x>$ and $Sym^{-1}(k(\overline{x}),\iota)=<\overline{x}>$.

The image of any one dimensional skew-hermitian form $<\alpha x>$ (which is automatically unramified) under the map $\partial_0$ is $<\overline{\alpha} \ \overline{x}>$. 

\textbf{Case B.2.2}: $\zeta=y$.

\textbf{Case B.2.2.1}: $\epsilon=1$.

We have $\pi'=z$. Hence $\sigma_{\pi'}=\sigma_z={\tau_{y}}_z=\tau_x$ is an orthogonal involution. 
Since $K(x)\subseteq Sym^1(D,\tau_y)=<1,x,z>$, the restriction of the involution $\tau_z$ on $K(x)$ is trivial. Therefore, $\overline{\sigma_{y}}=id$ so $W^1(k(\overline{x}),\overline{\tau_y})\cong W(k(x))$. Furthermore, since $K(x)\subseteq Sym^1(D,\tau_z)=<1,x,y>$, the restriction of the involution $\tau_z$ on $K(x)$ is trivial. Therefore, $\overline{\sigma_{\pi'}}=id$ so $W^1(k(\overline{x}),\overline{\tau_z})\cong W(k(x))$. In conclusion we have:

$$W^{1}((\frac{a,\pi}{K}),\tau_y) \cong  W(k(\overline{x})) \oplus  W(k(\overline{x})).$$

Note that $Sym^{1}((\frac{a,\pi}{K}),\tau_y)=<1,x,z>$. By Remark \ref{41} part (5), any one dimensional unramified hermitian form can be diagonalized as $<\delta +\alpha x>$ where $\delta + \alpha x$ is a unit of the field $K(x)$. The image of this form under $\partial_0$ is the quadratic form $<\overline{\delta} + \overline{\alpha} \ \overline{x}>$. Also by Remark \ref{43} part (3), any one dimensional ramified hermitian form can be diagonalized as $<\gamma z>$ where $\gamma \in {\mathcal{O}_K}^*$. The image of this form under $\partial_1$ is the quadratic form  $<\overline{\gamma} \ \overline{x}>$.

\textbf{Case B.2.2.2}: $\epsilon=-1$.

We have $\pi'=y$. Hence $\sigma_{\pi'}=\sigma_{y}=\tau_{yy}=\tau_{\pi}=\tau$ is the symplectic involution. The restriction of $\tau$ on 
$K(x)$ is the map $\iota$, and by Remark \ref{41} part (5), every skew-hermitian form over $(\frac{a,\pi}{K}, \tau_y)$ is ramified, hence we have the following isomorphism:

 $$W^{-1}((\frac{a,\pi}{K}),\tau_x) \cong W^{1}(k(\overline{x}),\iota).$$

Note that $Sym^{-1}((\frac{a,\pi}{K}),\tau_x)=<x>$ and $Sym^1(k(\overline{x}),\iota)=k$. By 
Remark \ref{43} part (4), any one dimensional ramified hermitian form can be diagonalized as $<\beta y>$, where $\beta \in {\mathcal{O}_K}^*$. The image of this form under $\partial_1$ is the hermitian form $<\overline{\beta}>$ over $(k(\overline{x}), \iota)$.\\ \

In the table below, we summarize  the results of this section: 

\begin{equation}\label{201}
    \begin{tabular}{||c  c c c c c c c c  c c ||} 
 \hline
 Case & $j$ & $\pi'$ & $\sigma$ & $\epsilon$ & $\pi''$ & $s_{\epsilon}$ &  $h_0$ & $h_1$  & $\overline{h_0}$ & $\overline{h_1}$   \\ [1.2ex] 
 \hline\hline
 A.1.1 & $1$ & $\pi$ & $\tau$   &  1 & $\pi$ & 1 &  $\delta$ & $\delta$ & $\overline{\delta}$ & $\overline{\pi^{-1}\delta}$   \\ 
 \hline
  A.1.2 & $1$ & $\pi x$ & $\tau$  &  -1 & $\pi x$ & $1$ &  $\alpha x + \beta y + \gamma z$ & $\alpha x + \beta y + \gamma z$  & $\overline{\alpha x + \beta y + \gamma z}$ & $\overline{\pi^{-1}(\alpha  -\gamma y + \beta z)}$   \\
 \hline
  A.2.1 & $1$ & $\pi$ & $\tau_x$  &  $1$ & $\pi$ & $1$ &  $\delta + \beta y + \gamma z$ & $\delta + \beta y + \gamma z$  & $\overline{\delta + \beta y + \gamma z}$ & $\overline{\pi^{-1}(\delta + \beta y + \gamma z)}$  \\
 \hline
  A.2.2 & $1$ & $\pi x$ & $\tau_x$ &  $-1$ & $\pi x$ & $1$ &  $\alpha x$ & $\alpha x$ & $\overline{\alpha} \ \overline{x}$ & $\overline{\pi^{-1}\alpha}$   \\
 \hline
  B.1.1 & $2$ & $y$ & $\tau$ & 
 $1$ & $\pi$ & $2$ &  $\delta$ & -  & $\overline{\delta}$ & -   \\
 \hline
  B.1.2 & $2$ & $y$ & $\tau$ &  $-1$ &  $y$ & $1$ &  $\alpha x$ & $\beta y + \gamma z$  & $\overline{\alpha} \ \overline{x}$ & $\overline{\beta} + \overline{\gamma} \ \overline{x}$  \\
 \hline
  B.2.1.1 & $2$ & $y$ & $\tau_x$ &  $1$ & $y$ & $1$ &  $\delta$ & $\beta y + \gamma z$ & $\overline{\delta}$ & $\overline{\beta} + \overline{\gamma} \ \overline{x}$   \\
 \hline
   B.2.1.2 & $2$ & $y$ & $\tau_x$ &  $-1$ & $\pi x$ & $2$ &  $\alpha x$ & - & $\overline{\alpha} \ \overline{x}$ & -    \\ \hline 
   B.2.2.1 & $2$ & $z$ & $\tau_y$ &  $1$ & $z$ & $1$ &  $\delta + \alpha x$ & $\gamma z$ & $\overline{\delta} + \overline{\alpha} \ \overline{x}$ & $\overline{\gamma} \ \overline{x}$  \\
 \hline
  B.2.2.2 & $2$  & $y$ & $\tau_y$ &  $-1$ & $y$ & $1$ & - & $\beta y$ & - & $\overline{\beta}$   \\ \hline

 \hline
\end{tabular}
\end{equation}

\end{document}